\documentclass[12pt,leqno,draft]{article} 
\usepackage{amsmath,amssymb,amsthm,amsfonts}
\usepackage{enumerate,color}
\makeatletter
\def\@cite#1#2{[{{\bfseries #1}\if@tempswa , #2\fi}]}
\renewcommand{\section}{%
\@startsection{section}{1}{\z@}
{0.5truecm plus -1ex minus -.2ex}%
{1.0ex plus .2ex}{\bfseries\large}}
\def\@seccntformat#1{\csname the#1\endcsname.\ }
\makeatother

\numberwithin{equation}{section} 
\newtheorem{thm}{Theorem}[section]
\newtheorem{corollary}[thm]{Corollary}
\newtheorem{lem}[thm]{Lemma}

\theoremstyle{definition}

\newtheorem{remark}{Remark}[section]

\newcommand{\ep}{\varepsilon}
\newcommand{\pa}{\partial}
\newcommand{\Rn}{\mathbb{R}^n}

\newcommand{\ol}[1]{\overline{#1}}

\newcommand{\tmax}{T_{\rm max}}
\newcommand{\lp}[2]{\Vert{#2}\Vert_{L^{#1}(\Omega)}}

\begin{document}
\footnote[0]
    {2010{\it Mathematics Subject Classification}\/. 
    Primary: 35K51; 
    Secondary: 35B45, 35A01, 92C17.
    }
\begin{center}
    \Large{{\bf 
           A unified method for boundedness 
           in fully parabolic chemotaxis systems 
           with signal-dependent sensitivity
          }}
\end{center}
\vspace{5pt}
\begin{center}
    Masaaki Mizukami\\
    \vspace{2pt}
    Department of Mathematics, 
    Tokyo University of Science\\
    1-3, Kagurazaka, Shinjuku-ku, Tokyo 162-8601, Japan\\
    {\tt masaaki.mizukami.math@gmail.com}\\
    \vspace{12pt}
    Tomomi Yokota%
   \footnote{Corresponding author}
   \footnote{Partially supported by Grant-in-Aid for
    Scientific Research (C), No.\,16K05182.}\\
    \vspace{2pt}
    Department of Mathematics, 
    Tokyo University of Science\\
    1-3, Kagurazaka, Shinjuku-ku, Tokyo 162-8601, Japan\\
    {\tt yokota@rs.kagu.tus.ac.jp}\\
    \vspace{2pt}
\end{center}
\begin{center}    
    \small \today
\end{center}
\vspace{2pt}
\newenvironment{summary}
{\vspace{.5\baselineskip}\begin{list}{}{%
     \setlength{\baselineskip}{0.85\baselineskip}
     \setlength{\topsep}{0pt}
     \setlength{\leftmargin}{12mm}
     \setlength{\rightmargin}{12mm}
     \setlength{\listparindent}{0mm}
     \setlength{\itemindent}{\listparindent}
     \setlength{\parsep}{0pt}
     \item\relax}}{\end{list}\vspace{.5\baselineskip}}
\begin{summary}
{\footnotesize {\bf Abstract.}
This paper deals with 
the Keller--Segel system 
\begin{equation*}
  \begin{cases}
    u_t=\Delta u - \nabla \cdot (u \chi(v)\nabla v), 
    & x\in\Omega,\ t>0, 
\\[1mm]	
    v_t=\Delta v + u - v, 
    & x\in\Omega,\ t>0, 
  \end{cases}
\end{equation*}
where $\Omega$ is a bounded domain in $\Rn$ 
with smooth boundary $\pa \Omega$, 
$n\geq 2$; $\chi$ is a function 
satisfying $\chi(s)\leq K(a+s)^{-k}$ 
for some $k\geq 1$ and $a\geq 0$. 
In the case that $k=1$, 
Fujie (J. Math.\ Anal.\ Appl.; 2015; 424; 675--684) 
established global existence of 
bounded solutions under the condition 
$K<\sqrt{\frac{2}{n}}$. 
On the other hand, when $k>1$, 
Winkler (Math.\ Nachr.; 2010; 283; 1664--1673) 
asserted global existence of bounded solutions 
for arbitrary $K>0$. 
However, there is a gap in the proof. 
Recently, Fujie tried modifying the proof; 
nevertheless it also has a gap. 
It seems to be difficult to show 
global existence of bounded solutions for arbitrary $K>0$. 
Moreover, the condition for $K$ when $k>1$ cannot 
connect to the condition when $k=1$. 
The purpose of the present paper is to 
obtain global existence and boundedness 
under more natural and proper 
condition for $\chi$ and 
to build a mathematical bridge between 
the cases $k = 1$ and $k >1$. 
}
\end{summary}
\newpage


\section{Introduction} 
%
%
The chemotaxis system proposed 
by Keller and Segel in \cite{K-S} 
describes a part of the life cycle of cellular 
slime molds with chemotaxis. 
After the pioneering work \cite{K-S}, 
a number of variations of the chemotaxis system 
are proposed and investigated (see e.g., 
\cite{B-B-T-W, H-P, H_until}). 

%
%
In this paper we consider the Keller--Segel system 
with signal-dependent sensitivity 
\begin{equation}\label{cp} 
  \begin{cases}
    u_t=\Delta u - \nabla \cdot (u \chi(v)\nabla v), 
    & x\in\Omega,\ t>0, 
\\[1mm]	
    v_t=\Delta v + u - v, 
    & x\in\Omega,\ t>0, 
\\[1mm]
    \nabla u\cdot \nu=\nabla v\cdot \nu = 0, 
    & x\in\pa \Omega,\ t>0,
\\[1mm]
    u(x,0)=u_0(x),\; v(x,0)=v_0(x), 
    & x\in\Omega,
  \end{cases}
\end{equation}
where $\Omega$ is a bounded domain in $\Rn$ 
($n\geq 2$) 
with smooth boundary $\pa \Omega$ 
and $\nu$ 
is 
the 
outward normal vector to $\pa\Omega$. 
The initial data $u_0$ and $v_0$ are
assumed to be nonnegative functions.
The unknown function $u(x,t)$ 
represents the population density of species and 
$v(x,t)$ shows the concentration of 
the 
substance 
at place $x$ and time $t$. 
As to the sensitivity function $\chi$, 
we are interested in 
functions 
generalizing 
%
\begin{align*}
  \chi(s)=\frac{K}{s}
\quad 
  \mbox{and} 
\quad
  \chi(s)=\frac{K}{(1+s)^2} 
\qquad
  (s>0). 
\end{align*}
In a mathematical view, 
the difficulty caused by the 
sensitivity function is 
to deal with the additional term $u\chi'(v)|\nabla v|^2$ 
which does not appear in the case that $\chi$ is a constant. 
Moreover, 
when $\chi$ is the singular sensitivity function, 
it is delicate to derive the estimate 
for $\frac{1}{v}$. 
In the case that 
$\chi(v)=\frac{K}{v}$, 
Winkler \cite{Winkler_2011} 
first attained 
global existence of classical solutions 
when $K<\sqrt{\frac{2}{n}}$ and 
global existence of weak solutions 
when $K<\sqrt{\frac{n+2}{3n-4}}$. 
However, the result in \cite{Winkler_2011} could not 
arrive at boundedness of solutions to \eqref{cp}. 
To overcome the difficulty in a singular sensitivity, 
Fujie \cite{Fujie_2015} established the 
uniform-in-time lower estimate for $v$, and 
show global existence of classical bounded solutions 
to \eqref{cp} in the case that 
$\chi(v)=\frac{K}{v}$ with
\begin{align}\label{condition;singular}
K<\sqrt{\frac{2}{n}}. 
\end{align}
As to the problem \eqref{cp} 
with $\chi(v)=\frac{K}{v}$ 
in the 2-dimensional setting, 
Lankeit \cite{Johannes_2016} obtained 
global existence of classical bounded solutions 
when $K<\chi_0$ ($\chi_0>1$). 
Futhermore Fujie--Senba \cite{F-S;p-p} dealt with 
the 2-dimensional problem 
which was replaced $v_t$ with $\tau v_t$ 
in \eqref{cp} and showed 
global existence and boundedness of radially symmetric solutions 
under the condition that $\chi(s)\to 0$ $(t\to\infty)$ 
and $\tau$ is sufficiently small. 
On the other hand, in the case that 
$\chi(v)\leq \frac{K}{(1+\alpha v)^k}$ 
$(k>1$, $\alpha>0$, $K > 0$), 
Winkler \cite{Winkler_2010ab} established 
global existence and boundedness in \eqref{cp} 
without any restriction on $K>0$.  
The result in \cite{Winkler_2010ab} 
affected the result in \cite{Fujie-Yokota} 
which mentioned global existence of classical bounded 
solutions to \eqref{cp} when $\chi$ is 
the strong singular sensitivity function 
such that 
$\chi(v)=\frac{K}{v^k}$ for $k>1$, 
and the method in \cite{Winkler_2010ab} 
was used in Zhang--Li \cite{Zhang-Li} 
and Zheng--Mu \cite{Zheng-Mu}. 
However, 
we cannot convince 
the result in \cite{Winkler_2010ab} 
that global existence and boundedness hold 
for all $K>0$, 
because there is a gap in the proof. 
Recently, Fujie tried modifying it;
nevertheless it also 
has a gap (cf.\ \cite{Fujie_DocTh}). 
In general, it seems to be difficult to show 
global existence of bounded solutions 
for arbitrary $K>0$. 
Moreover, if $k\to 1$, then 
the condition ``arbitrary $K>0$'' 
cannot connect to 
\eqref{condition;singular}. 

%
The purpose of this paper is 
to obtain global existence and boundedness 
in \eqref{cp} under a more natural and proper condition for $\chi$ 
and to build a mathematical bridge between 
the cases $k=1$ and $k>1$. 
We shall suppose that 
$\chi$ satisfies that 
\begin{align}\label{chiupper}
\chi(s)\leq \frac{K}{(a+s)^k}
\end{align}
with some $k\geq 1$, $a\geq 0$ and $K>0$ 
fulfiling  
\begin{align}\label{main} 
  K < k(a+\eta)^{k-1}\sqrt{\frac{2}{n}}. 
\end{align} 
Here 
\begin{align}\label{defeta}
  \eta:=
  \sup_{\tau>0}\left(
  \min\left\{
  e^{-2\tau} \min_{x\in\overline{\Omega}}v_0(x),\ 
  c_0\lp{1}{u_0}
  (1-e^{-\tau})
  \right\}
  \right),
\end{align} 
where $c_0>0$ is a lower bound for 
the fundamental solution of 
$w_t=\Delta w - w$ with Neumann boundary condition 
(for more detail, 
see Remark \ref{convex}). 
We suppose that 
\begin{align}\label{ini}
  0\leq u_0\in C(\ol{\Omega})\setminus \{0\} 
\quad 
  \mbox{and} 
\quad 
  \begin{cases}
  0< v_0\in W^{1,q}(\Omega)\ (\exists\, q>n)
&
  (a=0),
\\
  0\leq v_0\in W^{1.q}(\Omega)\setminus \{0\}
  \ (\exists\, q>n) 
&
  (a>0).
  \end{cases}
\end{align}
 
%
%
Now the main results read as follows. 
\begin{thm}\label{mainth}
 Let $n\geq 2$ and let $\Omega\subset \Rn$ be 
 a bounded domain with smooth boundary. 
 Assume that $\chi$ satisfies \eqref{chiupper} 
 with some $k\geq 1$, $a\geq 0$, $K>0$ fulfiling 
 \eqref{main}. 
 Then for any $u_0,v_0$ 
 satisfying \eqref{ini} with some $q>n$, 
 there exists an exactly one pair $(u,v)$ of 
 functions 
 \begin{align*}
   &u,\; v\in C(\ol{\Omega}\times [0,\infty))\cap 
   C^{2,1}(\ol{\Omega}\times (0,\infty))
 \end{align*} 
 which solves \eqref{cp}. 
 Moreover, the solution $(u,v)$ is uniformly bounded, 
 i.e., 
 there exists a constant $C>0$ such that 
 \begin{align*} 
   \lp{\infty}{u(\cdot,t)}+\lp{\infty}{v(\cdot,t)} 
   \leq C
\quad 
   \mbox{for all }t\geq 0. 
 \end{align*} 
 \end{thm}
\begin{remark} 
The unified condition 
\eqref{chiupper}
with $K>0$ satisfying \eqref{main} 
may be a natural condition for $\chi$. 
Indeed, when $k=1$, 
\eqref{main} 
becomes $\eqref{condition;singular}$: 
\begin{align*}
  K<\left.\left(k(a+\eta)^{k-1}\sqrt{\frac{2}{n}}\right)\right|_{k=1}= \sqrt{\frac{2}{n}}.
\end{align*}
\end{remark}
%
%
The main theorem tells us the result 
in the typical case of singular sensitivity.
\begin{corollary}
Let $\chi(s)=\frac{K}{s}$ with $K<\sqrt{\frac{2}{n}}$. 
Then for any $u_0,v_0$ 
satisfying \eqref{ini} with some $q>n$, 
\eqref{cp} has 
a unique global bounded classical solution. 
\end{corollary}

The strategy for the proof of Theorem \ref{mainth} 
is to construct the estimate for $\int_{\Omega}u^p$ 
with some $p>\frac{n}{2}$. 
One of the keys for this strategy is 
to derive the unified inequality 
\begin{align*}
  \frac{d}{dt}\int_{\Omega}u^p\varphi(v)
  &+\ep p(p-1)\int_\Omega u^{p-2}\varphi(v)|\nabla u|^2
  \leq 
  c\int_\Omega u^p\varphi(v)
  -r\int_\Omega u^{p+1}\frac{\varphi(v)}{(a+v)^k} 
\end{align*}
for some $\ep\in[0,1)$ and $c>0$, 
where 
  $$\varphi(s):=
  \exp\left\{-r\int_0^s\frac{1}{(a+s)^k}\,d\tau\right\}$$ 
with $r>0$. 
This function $\varphi$ 
constructed in \cite{M-Y_JDE, N-T_SIAM} 
unifies 
mathematical structures in 
the cases $k>1$ and $k=1$.

%
%
This paper is organized as follows. 
In Section 2 we collect basic facts 
which will be used later. 
In Section 3 
we give a unified view point in energy estimates. 
Section 4 is devoted to the 
proof of global existence and boundedness (Theorem 1.1).


\section{Preliminaries}
%
%
%
%
In this section 
we will collect elementary results. 
We first recall the uniform-in-time lower estimate for 
$v$ established by Fujie 
\cite{Fujie_2015, Fujie_DocTh}. 
\begin{lem}\label{lowesv} 
Let $u_0\in C(\ol{\Omega})$ 
and $u\in C(\ol{\Omega}\times [0,T))$ be 
nonnegative functions such that $u_0\not\equiv 0$ 
and $\int_\Omega u(\cdot,t)
=\int_\Omega u_0$ 
$(t\in[0,T))$. 
If  
$v_0\in C(\ol{\Omega})$ is a positive function 
in $\ol{\Omega}$ and 
$v\in C^{2,1}(\ol{\Omega}\times (0,T)) 
\cap C(\ol{\Omega}\times [0,T))$ 
is a classical solution of 
\begin{align*}
\begin{cases}
  v_t=\Delta v - v + u, 
&
  x\in\Omega,\ t\in (0,T),
\\
  \nabla v\cdot \nu =0 
&
  x\in\pa\Omega,\ t\in(0,T), 
\\
  v(x,0)=v_0(x)
&
  x\in \Omega, 
\end{cases}
\end{align*}
then for all $t\in(0,T)$, 
\begin{align*}
\inf_{x\in \Omega}v(x,t)\geq \eta>0, 
\end{align*}
where $\eta>0$ is defined as \eqref{defeta}. 
%
\end{lem}
\begin{remark}\label{convex}
When $\Omega$ is a convex bounded domain, 
the proof of this lemma is given in 
\cite[Lemma 2.2]{Fujie_2015} and 
the constant $\eta$ can be explicitly represented as 
\begin{align*}
  \eta:=
  \sup_{\tau>0}\left(
  \min\left\{
  e^{-\tau} \min_{x\in\overline{\Omega}}v_0(x),\ 
  \lp{1}{u_0}\cdot\int_0^\tau 
  \frac{1}{(4\pi r)^{\frac{n}{2}}}
  e^{-\left(r+\frac{\left(\max_{x,y\in\ol{\Omega}}|x-y|\right)^2}{4r}\right)}\,dr
  \right\}
  \right). 
\end{align*} 
On the other hand, if we do not 
assume the convexity of $\Omega$, 
then using the positivity of 
the fundamental solution $U(t,x;s,y)$ to 
$  w_t = \Delta w - w$ in $\Omega\times (0,T)$ 
with $\nabla w\cdot \nu = 0$ on $\pa \Omega \times (0,T)$ 
(see e.g., \cite{Ito}), 
we have that there exists $c_0>0$ such that 
for all $\tau \in (0,\frac{T}{2})$, 
\begin{align*}
  U(s+\tau,x; s,y)=U(\tau, x; 0, y)
  \geq c_0>0 
\quad
  \mbox{for all }x,y\in \ol{\Omega},\ s>0. 
\end{align*}
Then we can see the conclusion of Lemma \ref{lowesv} 
by a similar argument as in \cite[Lemma 2.2]{Fujie_2015}. 
\end{remark}
%
%
%
%
We next recall the well-known result about 
local existence of solutions to \eqref{cp} 
(see \cite[Theorem 2.3]{Johannes_2016}, 
\cite[Proposition 2.2]{F-S;p-p} and 
\cite[Lemma 2.1]{Winkler_2010ab}). 
\begin{lem}\label{localexistence}
Assume that $\chi$ satisfies \eqref{chiupper} and 
the initial data $u_0,v_0$ fulfil \eqref{ini} 
for some $q>n$. 
Then there exist  $\tmax\in (0,\infty]$ and exactly 
one pair $(u,v)$ of nonnegative functions 
\begin{align*}
  &u\in C(\ol{\Omega}\times [0,\tmax))
  \cap C^{2,1}(\ol{\Omega}\times(0,\tmax)),
\\
  &v\in C(\ol{\Omega}\times [0,\tmax))
  \cap C^{2,1}(\ol{\Omega}\times(0,\tmax))
  \cap L^\infty_{\rm loc}([0,\tmax);W^{1,q}(\Omega))
\end{align*} 
which solves \eqref{cp} in the classical sense 
and satisfies the mass conservation 
\begin{align}\label{masscon}
  \int_\Omega u(\cdot,t)=\int_\Omega u_0
\quad
  \mbox{for all }t\in (0,\tmax)
\end{align}
and the lower estimate 
\begin{align}\label{lowerestiv}
  \inf_{x\in \Omega}v(x,t)\geq \eta
\quad
  \mbox{for all }t\in (0,\tmax), 
\end{align}
where $\eta>0$ is defined as \eqref{defeta}. 
Moreover, if $\tmax<\infty$, then 
\begin{align*}
  \lim_{t\nearrow \tmax} 
  (\lp{\infty}{u(\cdot,t)} + 
  \|v(\cdot,t)\|_{W^{1,\infty}(\Omega)})=\infty. 
\end{align*}
\end{lem}
%
%
%
%
At the end of this section 
we shall recall the result about 
the estimate for $v$ 
in dependence on boundedness features of $u$ 
derived by a straightforward 
application of well-known smoothing 
estimates for the heat semigroup under homogeneous 
Neumann boundary conditions 
(see {\cite[Lemma 2.4]{Winkler_2011}} and 
\cite[Lemma 2.4]{Fujie_2015}). 
\begin{lem}\label{p-q}
 Let $T>0$ and\/ $1\leq \theta,\mu\leq \infty$. 
   If\/ $\frac{n}{2}(\frac{1}{\theta}-\frac{1}{\mu})<1$, 
   then there exists $C_1(\mu,\theta)>0$ such that 
   \begin{align*}
     \lp{\mu}{v(\cdot,t)}\leq C_1(\mu,\theta)
     \left(1+\sup_{s\in(0,\infty)}\lp{\theta}{u(\cdot,s)}\right)
   \quad
     \mbox{for all }t\in (0,T). 
   \end{align*}
%
\end{lem}
%
%
%

%
%

\section{A unified view point in energy estimates}
%
%
%
%
Let $(u,v)$ be the solution 
of \eqref{cp} on $[0,\tmax)$ 
as in Lemma \ref{localexistence}. 
For the proof of Theorem \ref{mainth} 
we will recall an useful fact to derive the 
$L^\infty$-estimate for $u$. 
\begin{lem}\label{Lp}
Assume that the solution $u$ of \eqref{cp} 
given in Lemma \ref{localexistence} satisfies 
\begin{align}\label{Lpesti}
  \lp{p}{u(\cdot,t)}\leq C(p) 
\quad 
  \mbox{for all }t\in (0,\tmax)
\end{align}
with some $p>\frac{n}{2}$ and $C(p)>0$. 
Then there exists a constant $C'>0$ such that 
\begin{align*}
  \lp{\infty}{u(\cdot,t)}\leq C' 
\quad 
  \mbox{for all }t\in (0,\tmax).
\end{align*}
\end{lem}
\begin{proof}
Combination of 
\eqref{lowerestiv} 
and the same argument as in 
\cite[Lemma 3.4]{Winkler_2011} 
leads to the conclusion of Lemma \ref{Lp}. 
\end{proof} 
%
%
%
%
\noindent {\bf Unified test function.} 
Thanks to Lemmas \ref{localexistence} 
and \ref{Lp} we will only make sure 
that the $L^p$-estimate for $u$ 
with some $p>\frac{n}{2}$ to 
show global existence and 
boundedness of solutions to \eqref{cp}. 
To establish \eqref{Lpesti} we introduce 
the functions $g$ and $\varphi$ by 
\begin{align}\label{defvarphi}
  g(s):=-r\int_\eta^s\frac{1}{(a+\tau)^k}\,d\tau, 
\quad 
  \varphi(s):=\exp\{g(s)\} 
\quad 
  (s\geq \eta), 
\end{align}
where 
$r>0$ is a constant fixed later 
and $\eta$ is defined as \eqref{defeta}. 
When $k>1$, by straightforward calculations 
we have 
\begin{align*}
  \varphi(s)=
  C_\varphi\exp\left\{\frac{r}{(k-1)(a+s)^{k-1}}\right\} 
\end{align*}
with $C_\varphi=\exp\{-r(k-1)^{-1}(a+\eta)^{-k+1}\}>0$, 
which is a similar function 
used in \cite{Winkler_2010ab}. 
On the other hand, when $k=1$, 
it follows that 
\begin{align*}
  \varphi(s)=\frac{C_\varphi}{(a+s)^{r}} 
\end{align*} 
with some constant $C_\varphi=(a+\eta)^r>0$, 
which is a similar function used 
in \cite{Fujie_2015}. 
Now we shall prove the following unified inequality 
by using the test function $\varphi(v)$. 
%
%
%
%
\begin{lem}\label{keylemma}
Assume that $\chi$ satisfies \eqref{chiupper}. 
Then for all $\ep\in [0,1)$ 
there exists $c>0$ such that 
\begin{align}\label{uniesti}
  \frac{d}{dt}\int_{\Omega}u^p\varphi(v)
  &+\ep p(p-1)\int_\Omega u^{p-2}\varphi(v)|\nabla u|^2
\\\notag
  &\leq 
  \int_{\Omega}u^pH_\ep(v)\varphi(v)|\nabla v|^2
  +c\int_\Omega u^p\varphi(v)
  -r\int_\Omega u^{p+1}\frac{\varphi(v)}{(a+v)^k}, 
\end{align}
where 
\begin{align*}
  H_\ep(s):= 
  \frac{(\ep p+1-\ep)r^2}{(1-\ep)(p-1)(a+s)^{2k}}
  +\frac{\ep pKr}{(1-\ep)(a+s)^{2k}}
  -\frac{kr}{(a+s)^{k+1}}
  +\frac{p(p-1)K^2}{4(1-\ep)(a+s)^{2k}}.
\end{align*}
\end{lem}
\begin{proof}
We let $p\geq 1$. 
From \eqref{cp} 
we have 
\begin{align}\label{firstineq;3.2}
  \frac{d}{dt}\int_\Omega u^p\varphi(v) 
  &=
  p\int_\Omega u^{p-1}\varphi(v)\nabla \cdot 
  (\nabla u-u\chi(v)\nabla v)
  +\int_\Omega u^p\varphi'(v)(\Delta v-v+u).
\end{align}
Integration by parts yields  
\begin{align}\label{cal;lem32}
  &p\int_\Omega u^{p-1}\varphi(v)\nabla \cdot 
  (\nabla u-u\chi(v)\nabla v)
  +\int_\Omega u^p\varphi'(v)\Delta v
\\\notag
  &= -p\int_\Omega \nabla(u^{p-1}\varphi(v))\cdot 
  (\nabla u-u\chi(v)\nabla v)
  -\int_\Omega \nabla(u^p\varphi'(v))\cdot \nabla v
\\\notag
  &=-p(p-1)\int_\Omega u^{p-2}
  \varphi(v)|\nabla u|^2
  +\int_\Omega u^{p-1}
  \left(p(p-1)\varphi(v)\chi(v)
  -2p\varphi'(v)\right)
  \nabla u\cdot\nabla v
\\\notag
  &\quad\ +\int_\Omega u^p (-\varphi''(v)+p\varphi'(v)\chi(v))|\nabla v|^2.
\end{align}
Invoking to the Young inequality, 
we infer that for all $\ep\in[0,1)$, 
\begin{align}\label{young;lemma3.2}
  &\int_\Omega u^{p-1}
  \left(p(p-1)\varphi(v)\chi(v)
  -2p\varphi'(v)\right)
  \nabla u\cdot\nabla v 
\\\notag
  &\leq (1-\ep)p(p-1)\int_\Omega 
  u^{p-2}\varphi(v)|\nabla u|^2
  +\int_\Omega u^p 
  \frac{(p(p-1)\varphi(v)\chi(v)-2p\varphi'(v))^2}
  {4(1-\ep)p(p-1)\varphi(v)}|\nabla v|^2. 
\end{align}
Thus combination of \eqref{firstineq;3.2}, 
\eqref{cal;lem32} and  \eqref{young;lemma3.2} 
yields that 
\begin{align}\label{step2;lem3.2}
  \frac{d}{dt}\int_\Omega u^p\varphi(v) 
  + \ep p(p-1)\int_\Omega 
  u^{p-2}\varphi(v)|\nabla u|^2
  \leq
  \int_\Omega u^{p}F_\varphi(v)|\nabla v|^2
  +\int_\Omega u^p\varphi'(v)(-v+u),
\end{align}
where 
\begin{align*}
  F_\varphi(s):=-\varphi''(s)
  -\frac{\ep}{1-\ep}p\varphi'(s)\chi(s)
  +\frac{p(p-1)}{4(1-\ep)}\chi(s)^2\varphi(s)
  +\frac{p\varphi'(s)^2}{(1-\ep)(p-1)\varphi(s)}. 
\end{align*}
Noting that 
\begin{align*}
  \varphi'(s)=g'(s)\varphi(s)
\quad
  \mbox{and}
\quad
  \varphi''(s)=g''(s)\varphi(s)+g'(s)^2\varphi(s),
\end{align*}
we can rewrite the function $F_\varphi(s)$ as
\begin{align*}
  F_\varphi(s)=\left(-g''(s)
  -\frac{\ep}{1-\ep}pg'(s)\chi(s)
  +\frac{p(p-1)}{4(1-\ep)}\chi(s)^2
  +\frac{(\ep p+1-\ep)g'(s)^2}{(1-\ep)(p-1)}
  \right)\varphi(s). 
\end{align*}
Recalling by \eqref{defvarphi} that 
\begin{align*}
  g'(s)=\frac{-r}{(a+s)^k}
\quad 
  \mbox{and}
\quad
  g''(s)=\frac{rk}{(a+s)^{k+1}},
\end{align*}
we see from \eqref{chiupper} that 
\begin{align}\label{Fineq}
F_\varphi(s)\leq H_\ep(s)\varphi(v),
\end{align}
where 
\begin{align*}
  H_\ep(s):= 
  \frac{(\ep p+1-\ep)r^2}{(1-\ep)(p-1)(a+s)^{2k}}
  +\frac{\ep pKr}{(1-\ep)(a+s)^{2k}}
  -\frac{kr}{(a+s)^{k+1}}
  +\frac{p(p-1)K^2}{4(1-\ep)(a+s)^{2k}}.
\end{align*}
Therefore we obtain from \eqref{step2;lem3.2} 
together with \eqref{Fineq} that 
\begin{align*}
  \frac{d}{dt}\int_{\Omega}u^p\varphi(v)
  &+\ep p(p-1)\int_\Omega u^{p-2}\varphi(v)|\nabla u|^2
\\
  &\leq 
  \int_{\Omega}u^pH_\ep(v)\varphi(v)|\nabla v|^2
  -r\int_\Omega u^{p}\varphi(v)\frac{(-v+u)}{(a+v)^k}. 
\end{align*}
We finally confirm from 
the boundedness of $\frac{s}{(a+s)^k}$ $(k\geq 1)$ 
that there exists $c>0$ satisfying 
\begin{align*}
  \int_\Omega u^p\varphi(v)\frac{rv}{(a+v)^k}
\leq 
  c\int_\Omega u^p\varphi(v), 
\end{align*}
and thus we obtain \eqref{uniesti}. 
\end{proof}
%
%
\section{Global existence and boundedness}
In this section we will show the $L^p$-estimate 
for $u$ with $p>\frac{n}{2}$ 
by using Lemma \ref{keylemma}. 
%
%
%
%
%
\subsection{Energy estimate in the case $k>1$}
We first derive the energy estimate in the case $k>1$. 
In this subsection we assume that $\chi$ satisfies 
\eqref{chiupper} with some $k>1$. 
Now we shall show 
the following inequality 
by modifying the method in \cite{M-Y_JDE}. 
%
%
%
%
%
%
\begin{lem}\label{energy;k>1}
Assume that \eqref{chiupper} 
and \eqref{main} are satisfied with some $k>1$, $a\geq 0$ and $K>0$. 
Then there exist 
$p>\frac{n}{2}$, $r>0$ and $\ep\in (0,1)$ such that 
\begin{align*} 
  H_\ep(s)\leq 0 
\quad
  \mbox{for all }s\geq \eta, 
\end{align*}
where $H_\ep$ is defined in Lemma \ref{keylemma} 
and  $\eta>0$ is defined as \eqref{defeta}, 
which implies that 
\begin{align}\label{step2;k>1}
  \frac{d}{dt}\int_{\Omega}u^p\varphi(v)
  &+\ep p(p-1)\int_\Omega u^{p-2}\varphi(v)|\nabla u|^2
  \leq 
  c\int_\Omega u^p\varphi(v)
  -r\int_\Omega u^{p+1}\frac{\varphi(v)}{(a+v)^k}. 
\end{align}
\end{lem}
\begin{proof}
We take $p\geq 1$, $r>0$ and $\ep\in[0,1)$ 
which will be fixed later. 
Due to the definition of $H_\ep$ we write as 
\begin{align*}
H_\ep(s)=a_1(s)r^2+a_2(s)r+a_3(s),
\end{align*}
where 
\begin{align*}
  &a_1(s):=
    \frac{\ep p+1-\ep}{(1-\ep)(p-1)(a+s)^{2k}}, 
\\
  &a_2(s):=
  \frac{\ep pK}{(1-\ep)(a+s)^{2k}} - 
  \frac{k}{(a+s)^{k+1}}, 
\\
  &a_3(s):=
    \frac{p(p-1)K^2}{4(1-\ep)(a+s)^{2k}}. 
\end{align*}
Noting from the condition \eqref{main} that 
there exist $p>\frac{n}{2}$ 
and $\ep\in(0,1)$ satisfying 
\begin{align*}
  K 
  <\frac{(1-\ep)k(a+\eta)^{k-1}}
  {\ep p+\sqrt{p(\ep p+1-\ep)}}, 
\end{align*}
we have that  
\begin{align*}
  \frac{\left(\ep p+\sqrt{p(\ep p+1-\ep)}\right)}{1-\ep}K
  \leq k(a+s)^{k-1} 
\quad 
  \mbox{for all }s\geq \eta. 
\end{align*}
%
This implies that 
the discriminant 
\begin{align*}
  D_r(s)&=a_2(s)^2-4a_1(s)a_3(s)
\\
  &=\left(
  \frac{\ep pK}{(1-\ep)(a+s)^{2k}} 
  - \frac{k}{(a+s)^{k+1}}\right)^2
  -\frac{p(\ep p+1-\ep)K^2}{(1-\ep)^2(a+s)^{4k}}
\\
  &=\frac{1}{(a+s)^{4k}}\left(
  \left(
  \frac{\ep pK}{(1-\ep)} 
  - k(a+s)^{k-1}\right)^2
  -\frac{p(\ep p+1-\ep)K^2}{(1-\ep)^2}
  \right)
\end{align*}
is nonnegative for all $s\geq \eta$. 
Finally, we show that 
there exists $r>0$ such that 
\begin{align}\label{aim;Hep<=0}
H_\ep(s)=a_1(s)r^2+a_2(s)r+a_3(s)\leq 0
\end{align}
for all $s\geq \eta$. 
Because $D_r(s)$ is nonnegative, we can define 
\begin{align*}
  r_{\pm}(s)&:= 
  \frac{-a_2\pm \sqrt{D_r(s)}}
  {2a_1(s)}
\\&=
  \frac{(1-\ep)(p-1)}{2(\ep p+1-\ep)}
  \left(
  -\frac{\ep pK}{(1-\ep)}+
  k(a+s)^{k-1}
  \pm (a+s)^{2k}
  \sqrt{D_r(s)}
  \right). 
\end{align*}
Then we see that $H_\ep(s)\leq 0$ 
for each $s\geq \eta$ and all 
$r\in(r_{-}(s),r_{+}(s))$. 
Since the functions 
\begin{align*}
  \widetilde{r}_{\pm}(\tau):=
  \tau \pm 
  \sqrt{\tau^2-\frac{p(\ep p+1-\ep)K^2}{(1-\ep)^2}}
\qquad
  \left(
  \tau \geq \frac{K\sqrt{p(\ep p+1-\ep)}}{1-\ep}
  \right)
\end{align*}
satisfy that
\begin{align*}
  \frac{d\widetilde{r}_+}{d\tau}(\tau)>0
\qquad
  \mbox{and}
\qquad
  \frac{d\widetilde{r}_-}{d\tau}(\tau)<0
\end{align*}
for all 
$\tau \geq \frac{K\sqrt{p(\ep p+1-\ep)}}{1-\ep}$, 
we obtain that 
\begin{align*}
  r_-(s)
\leq
  \widetilde{r}_-\left(
  \frac{K\sqrt{p(\ep p+1-\ep)}}{1-\ep}
  \right)
=
  \widetilde{r}_+\left(
  \frac{K\sqrt{p(\ep p+1-\ep)}}{1-\ep}
  \right)
\leq 
  r_+(s)
\end{align*}
holds for all $s\geq \eta$. 
Therefore if we put 
\begin{align*}
  r_0:=
    \widetilde{r}_-\left(
  \frac{K\sqrt{p(\ep p+1-\ep)}}{1-\ep}
  \right)
  =
  \frac{(p-1)K^2}{2}
  \sqrt{
    \frac{p}{\ep p+1-\ep}
  }
  >0,
\end{align*} 
then $r=r_0\in (r_-(s),r_+(s))$ for all $s\geq\eta$, 
which means that \eqref{aim;Hep<=0} holds for all $s\geq \eta$. 
This implies the end of the proof. 
\end{proof}
%
%
%
%
Now we are ready to show the $L^p$-estimate 
in the case $k>1$. 
\begin{lem}\label{Lp;k>1}
Assume that \eqref{chiupper} 
and \eqref{main} are satisfied with some $k>1$, $a\geq 0$ and $K>0$. 
Then there exist 
$p>\frac{n}{2}$ and $C>0$ such that 
\begin{align*}
  \lp{p}{u(\cdot,t)}\leq C
\quad
  \mbox{for all }t>0.
\end{align*}
\end{lem}
\begin{proof}
The proof is similar as in \cite{Winkler_2010ab}. 
From Lemmas \ref{keylemma} and \ref{energy;k>1} 
we obtain \eqref{step2;k>1} with some $p>\frac{n}{2}$, 
$r>0$ and $\ep\in (0,1)$. 
We shall show the $L^p$-estimate for $u$ 
by using \eqref{step2;k>1}. 
From the positivity of $u$, $v$ and $\varphi$ 
we have that 
\begin{align}\label{u^p+1;case1}
  -r\int_\Omega u^{p+1}\frac{\varphi(v)}{(a+v)^k}
  \leq 0.
\end{align}
We next deal with the term 
$\int_\Omega u^{p-1}\varphi(v)|\nabla u|^2$. 
Noting 
the boundedness of $\varphi$: 
\begin{align}\label{bddvarphi}
  C_\varphi
  \leq
  \varphi(s)
  \leq
  1
  \quad
  (s\geq\eta)
\end{align}
and combining the Gagliardo--Nirenberg inequality with 
the mass conservation property \eqref{masscon}, 
we deduce that there exists $c_1>0$ such that 
\begin{align}\label{G-N;step1}
  \int_\Omega u^p\varphi(v) 
  &\leq 
  \int_\Omega u^p 
  =
  \lp{2}{u^\frac{p}{2}}^2
\\\notag
  &\leq 
  c_1\left(\lp{2}{\nabla u^{\frac{p}{2}}}
  +\lp{\frac{2}{p}}{u^\frac{p}{2}}\right)^{2a}
  \lp{\frac{2}{p}}{u^\frac{p}{2}}^{2(1-a)}
\\\notag
  &=
  c_1\left(\lp{2}{\nabla u^{\frac{p}{2}}}
  +\lp{1}{u_0}^\frac{p}{2}\right)^{2a}
  \lp{1}{u_0}^{p(1-a)}
\end{align}
with $a=\frac{\frac{2n}{p}-
\frac{n}{2}}{\frac{2n}{p}+1-\frac{n}{2}}\in(0,1)$. 
From \eqref{step2;k>1}, \eqref{u^p+1;case1} 
and combination of \eqref{G-N;step1} and 
\begin{align*}
  C_\varphi\int_\Omega |\nabla u^\frac{p}{2}|^2 
\leq
  \int_\Omega u^{p-2}\varphi(v)|\nabla u|^2
\end{align*}
we see that there exist $c_2,c_3>0$ such that 
\begin{align*}
  \frac{d}{dt}\int_\Omega u^p\varphi(v) 
\leq 
  c \int_\Omega u^p\varphi(v)
  -c_2
  \left(\int_\Omega u^p\varphi(v)\right)^{\frac{1}{a}}
  +c_3,
\end{align*} 
which implies that there exist $p>\frac{n}{2}$, $r>0$ 
(determined in Lemma \ref{energy;k>1}) and $C>0$ satisfying  
\begin{align*}
  \int_\Omega u^p\varphi(v)\leq C. 
\end{align*}
Therefore we obtain from \eqref{bddvarphi} that 
\begin{align*}
\int_\Omega u^p \leq CC_\varphi^{-1}. 
\end{align*} 
Thus we attain the goal of the proof. 
\end{proof}
%
%
%
%
%
\subsection{Energy estimates in the case $k=1$} 
In this section we assume that 
the sensitivity function 
$\chi$ satisfies \eqref{chiupper} with $k=1$. 
We first show the following estimate for $H_0(s)$. 
\begin{lem}\label{energy;k=1}
Assume that 
\eqref{chiupper} and \eqref{main}
are satisfied with $k=1$ and with some $a\geq 0$ and $K>0$. 
Then for all $p\in (1,\frac{1}{K^2})$ 
there exists an interval $I_p$ such that 
for all $r\in I_p$ and $s\geq \eta$, 
\begin{align*}
  H_0(s)\leq 0, 
\end{align*}
where $H_0$ is defined as 
in Lemma \ref{keylemma} with $\ep=0$, 
which means 
\begin{align}\label{step2;k=1}
  \frac{d}{dt}\int_{\Omega}u^p\varphi(v)
  \leq 
  c\int_\Omega u^p\varphi(v)
  -r\int_\Omega u^{p+1}\frac{\varphi(v)}{a+v}. 
\end{align}
\end{lem}
\begin{proof}
We pick $p\geq 1$ and $r>0$ which will be fixed later. 
Due to the definition of $H_0(s)$ we write as
\begin{align*}
  H_0(s)&=
  \frac{1}{(a+s)^2}
  \left(\frac{1}{p-1}r^2- r 
  +\frac{p(p-1)K^2}{4}\right)
  \\
  &=:\frac{b_1r^2+b_2r+b_3}{(a+s)^2}. 
\end{align*}
From \eqref{main} we note that 
$(1,\frac{1}{K^2})$ is not an empty set. 
If we choose $p\in (1,\frac{1}{K^2})$, 
then the discriminant of $b_1r^2+b_2r+b_3$ 
is nonnegative: 
\begin{align*}
  D_r = 1-pK^2 \geq 0. 
\end{align*}
Thus we can define the interval $I_p$ as 
\begin{align*}
  I_p:=\left(\frac{p-1}{2}\left(1-\sqrt{1-pK^2}\right),
  \frac{p-1}{2}\left(1+\sqrt{1-pK^2}\right)\right)
\end{align*}
for each $p\in (1,\frac{1}{K^2})$. 
Then we have from straightforward calculations that 
for each $p\in (1,\frac{1}{K^2})$ and 
all $r\in I_p$, $b_1r^2+b_2r+b_3\leq 0$ holds, 
which yields that for all $p\in (1,\frac{1}{K^2})$ 
there exists the interval $I_p$ such that 
\begin{align*}
H_0(s)\leq 0 
\end{align*}
for all $r\in I_p$. 
This implies the end of the proof. 
\end{proof}
%
%
%
%
We next show the estimate for 
$\int_\Omega u^p\varphi(v)$. 
However we cannot easily obtain 
the estimate for $\int_\Omega u^p$ 
because $\varphi$ is not a bounded function. 
Therefore we will show the following lemma, 
which has an important 
role for obtaining the $L^p$-estimate. 
\begin{lem}\label{fromstep2}
Assume that 
\eqref{chiupper} and \eqref{main}
are satisfied with $k=1$ and with some $a\geq 0$ and $K>0$. 
Suppose that 
$p\in (1,\frac{1}{K^2})$, 
$r\in I_p$ such that $p-r\geq 1$. 
If there exists $C>0$ such that 
\begin{align}\label{p-r}
  \lp{p-r}{v(\cdot,t)}\leq C 
\quad
  \mbox{for all }t>0, 
\end{align}
then there exists $C_2(p,r)>0$ satisfying 
\begin{align*}
  \int_\Omega u^p\varphi(v)\leq C_2(p,r) 
\quad 
  \mbox{for all }t>0. 
\end{align*}
\end{lem}
\begin{proof} 
 The proof is similar as in \cite{Fujie_2015}. 
 We let $p\in (1,\frac{1}{K^2})$. 
 We denote by $I_p$ the interval 
 defined in Lemma \ref{energy;k=1}, 
 and choose $r\in I_p$. 
 We shall attain the conclusion 
 from \eqref{step2;k=1}. 
 By virtue of the H\"older inequality, 
 we infer that 
 \begin{align*}
   \int_\Omega u^p\varphi(v)
   &=
   C_\varphi^{-\frac{p}{p+1}}
   \int_\Omega 
   \left(
   u^{p+1}\frac{C_\varphi}{(a+v)^{r+1}}
   \right)^{\frac{p}{p+1}}\cdot
   (a+v)^{-r+\frac{p(r+1)}{p+1}}
 \\
   &\leq
   C_\varphi^{-\frac{p}{p+1}}
   \left(
   \int_\Omega u^{p+1}
   \frac{C_\varphi}{(a+v)^{r+1}}
   \right)^{\frac{p}{p+1}}
   \left(
   \int_\Omega (a+v)^{p-r}
   \right)^{\frac{1}{p+1}}
   \\
   &=
   C_\varphi^{-\frac{p}{p+1}}
   \left(
   \int_\Omega u^{p+1}
   \frac{\varphi(v)}{a+v}
   \right)^{\frac{p}{p+1}}
   \left(
   \int_\Omega (a+v)^{p-r}
   \right)^{\frac{1}{p+1}}.
 \end{align*}
Noting from \eqref{p-r} 
and the fact $p-r\geq 1$ that 
\begin{align*}
  \lp{p-r}{a+v(\cdot,t)}\leq a+C 
  \quad 
  \mbox{for all }t>0, 
\end{align*} 
we obtain that there exists $c_1>0$ such that 
\begin{align}\label{p+1;k=1}
  \int_\Omega u^p\varphi(v)
\leq 
  c_1 \left(
  \int_\Omega u^{p+1}
  \frac{\varphi(v)}{a+v}\right)^{\frac{p}{p+1}}. 
\end{align}
Plugging \eqref{p+1;k=1} into \eqref{step2;k=1} 
yields that there exists $c_2>0$ such that 
\begin{align*}
  \frac{d}{dt}\int_{\Omega}u^p\varphi(v)
  \leq 
  c\int_\Omega u^p\varphi(v)
  -c_2\left(
  \int_\Omega u^{p}\varphi(v)\right)^{\frac{p+1}{p}}. 
\end{align*}
Therefore we have from a standard ODE comparison 
argument that there exists $C>0$ such that 
\begin{align*}
\int_\Omega u^p\varphi(v)\leq C,
\end{align*}
which means the end of the proof. 
\end{proof}
%
%
%
%
\begin{lem}\label{Lp;k=1}
Assume that 
\eqref{chiupper} and \eqref{main}
are satisfied with $k=1$ and with some $a\geq 0$ and $K>0$. 
Then there exist 
$p>\frac{n}{2}$ and $C>0$ such that 
\begin{align*}
  \lp{p}{u(\cdot,t)}\leq C.
\end{align*}
\end{lem}
\begin{proof}
We use a similar argument as in \cite[Proof of Main Theorem (Step 1)]{Fujie_2015}. 
First we choose a pair $(p_0,r_0)$ such that 
\begin{equation*}
\left\{
\begin{array}{l}
  p_0\in
  \left(
  1,
  \min\left\{\dfrac{1}{K^2}, 
  n+1, 
  \dfrac{n+2}{n-2}\right\}
  \right),\\
  r_0:=\dfrac{p_0-1}{2}.
\end{array}
\right.
\end{equation*}
Then $p_0>r_0$, $r_0<\frac{n}{2}$, $r_0\in I_p$, 
$p_0-r_0\geq 1$ and $p_0-r_0<\frac{n}{n-2}$. 
Since $\frac{n}{2}(1-\frac{1}{p_0-r_0})<1$ from 
the inequality $p_0-r_0<\frac{n}{n-2}$, 
it follows from 
Lemma \ref{p-q} and \eqref{masscon} that there exists 
$b_0>0$ such that 
\begin{align*}
  \lp{p_0-r_0}{v(\cdot,t)}
  \leq C_1(p_0-r_0,1)
  \left(1+\sup_{s\in (0,\infty)}\lp{1}{u(\cdot,s)}\right)
  \leq b_0.
\end{align*}
Thus we have from the fact $p_0-r_0\geq 1$ and 
Lemma \ref{fromstep2} that 
\begin{align*}
  \int_\Omega u^{p_0}(a+v)^{-r_0}
  \leq   C_\varphi^{r_0} C_2(p_0,r_0). 
\end{align*}
We will show that for all 
$q_0\in (\frac{2p_0}{p_0+1},
\min\{p_0,\frac{n(p_0-r_0)}{n-2r_0}\})$ 
there exists $b_0'>0$ such that 
\begin{align}\label{q_0esti}
\int_\Omega u^{q_0}\leq b_0'. 
\end{align}
We first note from the fact $p_0\not\in (\frac{n}{n+4},1)$ 
that 
\begin{align*}
  \frac{2p_0}{p_0+1}\leq p_0 
\quad 
  \mbox{and}
\quad
  \frac{2p_0}{p_0+1} \leq \frac{n(p_0+1)}{2(n-p_0+1)}
  =\frac{n(p_0-r_0)}{n-2r_0}. 
\end{align*}
Now we fix $q_0\in (\frac{2p_0}{p_0+1},
\min\{p_0,\frac{n(p_0-r_0)}{n-2r_0}\})$. 
Applying the H\"older inequality, we infer that 
\begin{align}\label{q_0holder}
  \int_\Omega u^{q_0}
  &\leq 
  \left(\int_\Omega u^{p_0}
  (a+v)^{-r_0}\right)^{\frac{r_0}{p_0}}
  \cdot \left(
  \int_\Omega (a+v)^\frac{q_0r_0}{p_0-q_0}
  \right)^\frac{p_0-r_0}{p_0}
\\\notag
  &\leq
  \left(
  C_\varphi^{r_0}C_2(p_0,r_0)
  \right)^\frac{r_0}{p_0}
  \cdot \lp{\frac{q_0r_0}{p_0-r_0}}{a+v(\cdot,t)}^\frac{q_0r_0}{p_0}. 
\end{align} 
We can confirm from the fact 
$\frac{n}{2}(\frac{1}{q_0}-\frac{p_0-q_0}{q_0r_0})<1$
and Lemma \ref{p-q} 
that
\begin{align}\label{q_0r_0}
\lp{\frac{q_0r_0}{p_0-r_0}}{v(\cdot,t)}
\leq C_1\left(\frac{q_0r_0}{p_0-r_0},q_0\right)
\left(1+\sup_{s\in (0,\infty)}\lp{q_0}{u(\cdot,s)}\right). 
\end{align} 
Therefore plugging \eqref{q_0r_0} and 
the inequality 
$\frac{q_0r_0}{p_0-r_0}\geq 1$ 
(from the fact $q_0\geq \frac{2p_0}{p_0+1}$) 
into \eqref{q_0holder} derives that 
there exists $b_0''>0$ satisfying 
\begin{align*}
  \lp{q_0}{u(\cdot,t)}
  \leq 
  b_0'' 
  \left(1 + 
  \sup_{s\in (0,\infty)}\lp{q_0}{u(\cdot,s)}\right)^{\frac{r_0}{p_0}}. 
\end{align*}
Noting that $\frac{r_0}{p_0}<1$, 
we obtain \eqref{q_0esti}. 
In the above argument, 
since $\frac{n}{2}<\min\{\frac{1}{K^2},n+1\}$, 
if $\frac{n+2}{n-2}>\frac{n}{2}$, 
then we can fix $p_0>\frac{n}{2}$. 
Thus we can find $q_0>\frac{n}{2}$ due to 
$p_0< \frac{n(p_0+1)}{2(n-p_0+1)}
=\frac{n(p_0-r_0)}{n-2r_0}$ 
if $\frac{n+2}{n-2}>\frac{n}{2}$. 
On the other hand, if $\frac{n+2}{n-2}\leq\frac{n}{2}$, 
then we proceed the iteration argument. 
For all $i\in \mathbb{N}$, 
we fix a pair $(p_i,r_i)$ defined 
inductively such as 
\begin{align*}
\left\{
\begin{array}{l}
  p_i\in \left(p_{i-1},
  \min\left\{\dfrac{1}{K^2},n+1,
  \dfrac{p_{i-1}(n+2)}{n-2p_{i-1}}\right\} \right),
\\
  r_i:=\dfrac{p_i-1}{2},
\end{array}
\right.
\end{align*}
then we can see that $p_i>r_i$, $r_i<\frac{n}{2}$ 
and $r_i\in I_{p_i}$. 
Moreover, 
since the relation  
$p_i<\frac{p_{i-1}(n+2)}{n-2p_{i-1}}$ 
implies that 
\begin{align*}
p_i-r_i < 
\frac
  {n\cdot \frac{n(p_{i-1}+1)}{2(n-p_{i-1}+1)}}
  {n-2\cdot\frac{n(p_{i-1}+1)}{2(n-p_{i-1}+1)}},
\end{align*}
we can find some $q_{i-1}\in (\frac{2p_{i-1}}{p_{i-1}+1},
\min\{p_{i-1},\frac{n(p_{i-1}+1)}{2(n-p_{i-1}+1)}\})$ 
satisfying 
\begin{align*}
p_i-r_i< \frac{nq_{i-1}}{n-2q_{i-1}}, 
\end{align*}
which means that 
$\frac{n}{2}(\frac{1}{q_{i-1}}-\frac{1}{p_i-r_i})<1$. 
Therefore we obtain 
\begin{align}\label{pi-ri}
  \lp{p_i-r_i}{v(\cdot,t)} 
  \leq  
  C_1(p_i-r_i,q_{i-1})
  \left(
  1+\sup_{s\in(0,\infty)}\lp{q_{i-1}}{u(\cdot,s)}
  \right).
\end{align}
According to \eqref{pi-ri} and $p_i-r_i\geq 1$ 
into Lemma \ref{fromstep2}, 
we have 
\begin{align*}
  \int_\Omega u^{p_i}(a+v)^{-r_i}
  \leq C_\varphi^{r_i} C_2(p_i,r_i). 
\end{align*}
Due to a similar argument to the first iteration, 
we can show that 
\begin{align*}
  \lp{q_i}{u(\cdot,t)}\leq b_i'(q_i)
\qquad 
  \mbox{for all } q_i\in \left(\frac{2p_{i}}{p_{i}+1},
\min\left\{p_{i},
\frac{n(p_{i}-r_{i})}{n-2r_{i}}\right\}\right) 
\end{align*}
with some constant $b_i'(q_i)>0$. 
Because the increasing function 
$f(x):=\frac{x(n+2)}{n-2x}$ satisfies 
$f(x)>1$ and $f(x)\to \infty$ as $x\to \frac{n}{2}$, 
we can find some $i_0\in\mathbb{N}$ such that 
$p_{i_0} > \frac{n}{2}$ and 
$p_{i_0-1}\leq \frac{n}{2}$, 
and hence $q_{i_0}>\frac{n}{2}$. 
Therefore we verify 
\begin{align*}
\lp{p}{u(\cdot,t)}\leq C
\end{align*}
with some $p>\frac{n}{2}$ and some $C>0$, 
which completes the proof. 
\end{proof}
%
%
%
\subsection{Proof of Theorem \ref{mainth}}
Combination of the $L^p$-estimate for $u$ 
(see Lemma \ref{Lp;k>1} or Lemma \ref{Lp;k=1}) 
and Lemma \ref{Lp} directly leads to 
Theorem \ref{mainth}. \qed


\newpage
\begin{small}

\end{small}

\end{document}